\documentclass[12pt]{article}
\usepackage{bez123,calc,curves,ebezier,epic,eepic,graphicx,multiply,rotating}
\everymath{\displaystyle}
\usepackage{graphicx}
\usepackage{pst-all}
\usepackage{color}
\usepackage{amssymb}
\usepackage{amsmath}
\newtheorem{prethm}{{\bf Theorem}}

\newenvironment{thm}{\begin{prethm}{\hspace{-0.5
               em}{\bf .}}}{\end{prethm}}
\newtheorem{prelemma}{{\bf Lemma}}

\newtheorem{preex}{{\bf Example}}

\newtheorem{preprop}{{\bf Proposition}}

\newenvironment{prop}{\begin{preprop}{\hspace{-0.5em}{\bf .}}}{\end{preprop}}
\newtheorem{precor}{{\bf Corollary}}

\newenvironment{cor}{\begin{precor}{\hspace{-0.5
               em}{\bf .}}}{\end{precor}}
\newtheorem{preremark}{{\bf Remark}}

\newenvironment{remark}{\begin{preremark}{\hspace{-0.5
               em}{\bf.}}}{\end{preremark}}
\newtheorem{preprob}{{\bf Problem}}

\newtheorem{predefin}{{\bf Definition}}

\newtheorem{preconj}{{\bf Conjecture}}

\newtheorem{preprobb}{{\bf Problem}}

\newtheorem{prelem}{{\bf Theorem}}

\newenvironment{proof}{{\bf Proof.}\rm }{\hfill{$\Box$}}

\newtheorem{presolution}{{\bf Solution.}}

\def\newpic#1{}
\def\qed{\ifhmode\unskip\nobreak\fi\quad\ifmmode\Box\else$\Box$\fi}
\title{\vspace{-0.1cm}\Large\bf On dynamic monopolies of graphs with general thresholds}
\author{\large\bf Manouchehr Zaker\footnote{E-mail: mzaker@iasbs.ac.ir}
\vspace{5mm}\\
    Department of Mathematics,\\
     Institute for Advanced Studies in Basic Sciences,\\
     Zanjan 45137-66731, Iran}

    \date{}

\begin{document}
\maketitle
\begin{abstract}
\noindent Let $G$ be a graph and ${\mathcal{\tau}}: V(G)\rightarrow \Bbb{N}$ be an assignment of thresholds to the vertices of $G$. A subset of vertices $D$ is said to be dynamic monopoly (or simply dynamo) if the vertices of $G$ can be partitioned into subsets $D_0, D_1, \ldots, D_k$ such that $D_0=D$ and for any $i=1, \ldots, k-1$ each vertex $v$ in $D_{i+1}$ has at least $t(v)$ neighbors in $D_0\cup \ldots \cup D_i$. Dynamic monopolies are in fact modeling the irreversible spread of influence such as disease or belief in social networks. We denote the smallest size of any dynamic monopoly of $G$, with a given threshold assignment, by $dyn(G)$. In this paper we first define the concept of a resistant subgraph and show its relationship with dynamic monopolies. Then we obtain some lower and upper bounds for the smallest size of dynamic monopolies in graphs with different types of thresholds. Next we introduce
dynamo-unbounded families of graphs and prove some related results. We also define the concept of a homogenious society that is a graph with probabilistic thresholds satisfying some conditions and obtain a bound for the smallest size of its dynamos. Finally we consider dynamic monopoly of line graphs and obtain some bounds for their sizes and determine the exact values in some special cases.
\end{abstract}



%%%%%%%%%%%%%%%%%%%%%%%%%%%%%%%%%%%%%%%%%%%%%%%%%%%%%%%%%%%%%%%%%%%%%%%%%%%

\section{Motivation and related works}

\noindent In recent years, great attentions have been paid to the modeling and analysis of the spread of belief or influence in complex networks. Various problems in social and virtual networks such as world wide web or models of distributed computing can be formalized in terms of the spread of influence. Elections in societies where individuals
decide whether to vote a certain candidate, spread of disease among people or virus in
world wide web or any web of computers are some examples of these problems. A network in all of these examples which is simply consisted of a set of elements (e.g. agents in
social networks or computing units in distributed computing systems) and some relationships or interactions between these elements can be conveniently modeled by a graph whose nodes represent the elements and edges represent the links of the network. For the graph theoretical notions, not defined in this paper, we refer the reader to \cite{BM}.

\noindent The model to be studied in this paper is as follows. A graph $G$ on the vertex set $V(G)$ and the edge set $E(G)$ together with an assignment of thresholds $\tau: V(G)\rightarrow \Bbb{N}$ to the vertices of $G$ is given. The discrete time dynamic process corresponding to the
threshold assignment $\tau$ is defined as follows:

\noindent The process starts with a subset $D$ of vertices which consists of the vertices having the state $+$ at time $0$. We denote the set of vertices of state $+$ in time $i$ by $D_i$. So at the beginning i.e. at time $0$ we have $D_0=D$. Then at any time $i+1\geq 1$, the state of any vertex $v$ changes to the state $+$ provided that at least $t(v)$ neighbors of $v$ belong to $D_i$. Also if the state of $v$ is already $+$ in time $i$ (i.e. when $v\in D_i$) then its state remains $+$ as before. If at a certain time $t$ of the process a vertex $v$ has state + then $v$ is said to be an active vertex. Note that the process defined above is progressive or irreversible i.e. when the state of a vertex becomes $+$ at some step of the process then its state remains unchanged until the end of the process.

\noindent By a $\tau$-dynamic monopoly we mean any subset $D$ of the vertices of $G$ such that by starting from $D$, all the vertices of $G$ get state $+$ at the end of the process.
Throughout the paper by $(G,\tau)$ we mean a graph $G$ together with a threshold assignment $\tau$ to the vertices of $G$. We simply write dynamic monopoly or (following some authors) dynamo instead of $\tau$-dynamic monopoly. By the size of a dynamo $D$ we mean the cardinality of $D$. It is easy to see that a subset of vertices $D$ in a graph $(G,\tau)$ is dynamo if and only if there exists a partition of $V(G)$ into subsets $D_0, D_1, \ldots, D_k$ such that $D_0=D$ and for any $i=1, \ldots, k-1$ each vertex $v$ in $D_{i+1}$ has at least $t(v)$ neighbors in $D_0\cup \ldots \cup D_i$. Dynamic monopolies have been widely studied by various authors. Some related graph theoretical and algorithmic results concerning dynamic monopolies have been obtained in \cite{DR,FKRRS}. Dynamic monopolies in terms of repetitive polling games were studied in \cite{P1}. More polynomial time or approximation algorithms were obtained in \cite{C}. Dynamic monopolies of special families of graphs were studied in \cite{FGS,FLLPS,LPS}. Also \cite{KKT} studies optimization formulations of dynamic monopolies and obtains some approximation algorithms. In \cite{P1}, controlling monopolies, a concept similar to dynamic monopolies has been introduced. The survey paper \cite{P2} surveys some of the results concerning various kinds of monopolies including dynamic monopolies. It also presents various applications of
these concepts in distributed computing and probabilistic polling models.

\noindent {\bf The outline of the paper is as follows.} In Section 2 we first introduce resistant subgraphs and show their relationships with dynamic monopolies. Then we obtain some lower and upper bounds for the dynamo size of graphs with various types of thresholds. Finally in Section 2 we determine the smallest dynamo of Generalized Petersen graph $GP(n,k)$. In Section 3 we introduce dynamo-unbounded families of graphs and obtain some results in this regard. We show the importance of the concept of dynamo-unbounded families by its applications in some famous social network problems. Homogenious societies are also defined in Section 3 and a result concerning their dynamo sizes is given in Section 3. In the last section we consider line graphs and
obtain some lower bounds for their dynamo numbers and determine the exact value for the line graph of the complete bipartite graphs.

\section{Some bounds for dynamo size of graphs}

\noindent We begin with the following concept. Given $(G,\tau)$, by a resistant subgraph of $G$ we mean any subgraph $K$ such that for any vertex
$v\in K$ one has $d_K(v)\geq d_G(v)-t(v)+1$, where $d_G(v)$ is the degree of $v$ in $G$. The following proposition provides a sufficient and necessary condition for graphs containing no resistant subgraphs.

\begin{prop}
A subgraph $H$ of $(G,\tau)$ does not contain any resistant subgraph of $G$ if and only if the vertices of $H$ can be labeled as $v_1, v_2, \cdots, v_n$ such that $v_i$ has at most $d_G(v_i)-t(v_i)$ neighbors among $\{v_i, v_{i+1}, \ldots, v_n\}$.\label{prop1}
\end{prop}

\noindent \begin{proof}
Assume first that $H$ contains no resistant subgraph. Then $H$ itself is not a resistant graph and so there exists a vertex $v_1\in H$ such that $d_H(v_1)\leq d_G(v_1)-t(v_1)$.
Set $H_1=H\setminus \{v_1\}$. Since $H_1$ too is not resistant then there exists $v_2$ such that $d_{H_1}(v_2)\leq d_G(v_2)-t(v_2)$. This means that $v_2$ has at most $d_G(v_2)-t(v_2)$ neighbors in $H[v_2, \ldots, v_n]$. We repeat this argument and obtain $v_1, v_2, \ldots, v_n$ such that $d_{H_{i-1}}(v_i)\leq d_G(v_i)-t(v_i)$ where $H_{i-1}=H[v_i, \ldots, v_n]$.

\noindent Assume now that the vertices of $H$ are labeled as specified in the proposition. Let by contrary that $H$ has a resistant subgraph $K$. Let $i$ be the smallest index with $v_i\in K$. Then $K\subseteq H_{i-1}=H[v_i, \ldots, v_n]$. This implies $d_K(v_i)\leq d_{H_{i-1}}(v_i)$. From one side we have $d_G(v_i)-t(v_i)+1\leq d_K(v_i)$ and from other side $d_{H_{i-1}}(v_i)\leq d_G(v_i)-t(v_i)$. This contradiction completes the proof.
\end{proof}

\begin{prop}
A subset $D$ in $(G,\tau)$ is dynamic monopoly if and only if $G\setminus D$ does not contain any resistant subgraph.\label{prop2}
\end{prop}

\noindent \begin{proof}
Assume first that there exists a subgraph $K$ of $G\setminus D$ which is resistant. Let $D_i$ be the set of vertices of $G$ which become active at time $i$, by starting from $D_0=D$. Let $v$ be any arbitrary vertex of $K$. Then $v$ has at most $t(v)-1$ neighbors in $G\setminus K$. Since $D\subseteq G\setminus K$ then this shows that $v$ has at most $t(v)-1$ neighbors in $D$. Therefore $K\cap D_1=\varnothing$ and in particular $v$ has at most $t(v)-1$ neighbors in
$D_1$. We repeat the argument we made above, for $v$ and $D_1$ and also for other $D_i$, $i=2, \ldots$ to conclude that $K$ remains outside $D_i$ for any $i$. This shows that $D$ can not be a dynamic monopoly.

\noindent Now assume that $H=G\setminus D$ does not contain any resistant subgraph. This shows that
$H$ itself is not resistant. Then there exists a vertex $v_1$ of $H$ with $d_H(v_1)\leq d_G(v_1)-t(v_1)$.
Namely $v_1$ has at least $t(v_1)$ neighbors in $K$. Then $v_1\in D_1$. Now we repeat this argument for $H\setminus v_1$ which is not a resistant subgraph of $H$ by the hypothesis on $H$. We obtain another vertex $v_2$ of $H$ with $v_2\in D_2$. Naturally this technique can be applied until all the vertices of $H$ get state $+$ in at most $|H|$ time steps. This completes the proof.
\end{proof}

\noindent The following useful remark comes immediately.

\begin{remark}
If a vertex $v\in G$ is such that $t(v)\geq d_G(v)+1$ then the subgraph of $G$ consisting of the single vertex $v$ is a resistant subgraph.
\end{remark}

\noindent Throughout the paper for any subset $S\subseteq V(G)$ we denote the subgraph of $G$ induced on $S$ by $G[S]$.

\begin{thm}
Let $D$ be a dynamic monopoly of size $k$ in $G$. Set $H=G\setminus D$ and let $t_{max}$ be the maximum threshold among the vertices of $H$. Then

$(i)~~  \sum_{v\in H}t(v) \leq |E(G)|-|E(G[D])|-\delta(G)+t_{max}$

$(ii)$~~ $\sum_{v\in H}t(v) \leq |E(G)|$ provided that $t(v)\leq d_G(v)$ for any vertex $v\in H$. \label{bound1}
\end{thm}

\noindent \begin{proof}
By Proposition \ref{prop2}, $H$ does not contain any resistant subgraph and so by Proposition \ref{prop1} the vertices of $H$ can be labeled as $v_1, \ldots, v_n$ in such a way that for any $i$, $d_{H_i}(v_i)\leq d_G(v_i)-t(v_i)$ where $H_i=H[v_i, \ldots, v_n]$.
In particular $H=H_1$ and $H_2=H\setminus v_1$. For the number of edges of $H_i$ we have
$$|E(H_i)| = |E(H_{i-1})|-d_{H_{i-1}}(v_{i-1})\geq |E(H_{i-1})| - d_G(v_{i-1}) + t(v_{i-1}).$$
\noindent We use recursively the above inequality and obtain the following
$$|E(H_i)|\geq |E(H)|-\sum_{j=1}^{i-1} d_G(v_j) + \sum_{j=1}^{i-1} t(v_j).$$
\noindent For $i=|H|$ we have $|E(H_i)|=0$ and obtain the following
$$|E(H)|\leq \sum_{v\in H} d_G(v) - d_G(v_n) - \sum_{v\in H} t(v) + t(v_n).~~~~~{\bf (1)}$$

\noindent Now we estimate the number of edges of $G$. Let $e$ be the number of edges between $D$ and $H$. We have $|E(G)|=|E(G[D])|+ e +|E(H)|$. Using {\bf (1)} we obtain
$$|E(G)|\leq |E(G[D])|+ e + \sum_{v\in H} d_G(v) - d_G(v_n) - \sum_{v\in H} t(v) + t(v_n).$$
\noindent We note that $\sum_{v\in H} d_G(v)=e+2|E(H)|$. Let $t_{max}$ be the maximum threshold among the vertices of $H$. We have now
$$|E(G)|\leq |E(G[D])|+ e +e +2|E(H)|-\delta(G)+t_{max}-\sum_{v\in H} t(v).$$
\noindent Therefore $$\sum_{v\in H} t(v)\leq |E(G)|-|E(G[D])|-\delta(G)+t_{max}.$$
\noindent To prove part (ii), note that by {\bf (1)} we have $|E(H)|\leq \sum_{v\in H} d_G(v) - \sum_{v\in H} t(v)$. Now by $|E(G)|=|E(G[D])|+ e +|E(H)|$ and $\sum_{v\in H} d_G(v)=e+2|E(H)|$ the desired inequality is obtained.
\end{proof}

\noindent As a corollary of Theorem \ref{bound1}, part (ii) we have the following result for regular graphs.

\begin{cor}
Let $G$ be a regular graph of degree $2r+1$. Let also $t(v)=r+1$ for any vertex of $G$. Then any dynamo for $G$ has at least $\frac{n+2r}{2(r+1)}$ vertices.\label{regular}
\end{cor}

\noindent In Theorem \ref{GP} we show that the bound obtained in Corollary \ref{regular}
is a tight bound when $r=1$ i.e. for cubic graphs. Of course we believe that it is also
tight for regular graphs of higher degrees. Before stating Theorem \ref{GP}, we present an upper bound for the dynamic monopoly of graphs in which the threshold of any vertex is 2.

\begin{thm}
Let $G$ be a graph on $n$ vertices such that no connected component of $G$ is isomorphic to an odd cycle, with $\delta(G)\geq 2$ and $t(v)=2$ for any $v\in G$. Let also $S$ be any domination set of $G$ and let $c$ be the number of connected components of $G\setminus S$.
Then the following bounds hold.

(i) $dyn(G)\leq n/2$

(ii) $dyn(G)\leq dyn(G[S])+c$

(iii) $dyn(G)\leq \frac{|S|}{2}+c$.\label{t=2}
\end{thm}

\noindent \begin{proof}
We may assume that $G$ is connected. To prove (i), let $\{C_1, \ldots, C_k\}$ be a set of vertex disjoint cycles in $G$ with the maximum
cardinality. This set is non-empty because $\delta(G)\geq 2$. In the following we obtain a dynamo denoted by $D$ of cardinality at most $n/2$. The subgraph $G\setminus (C_1\cup \ldots \cup C_k)$ is a forest, denote it by $F$. Let $T$ be any connected component of $F$. The tree $T$ contains a domination set $S$ of the cardinality at most $|T|/2$. From any connected component $T$ of $F$ we consider a minimum domination, say $S$ for $T$ and add the vertices of $S$ to $D$. Note that $|S|\leq |T|/2$.

\noindent Now we consider the odd cycles of $\{C_1, \ldots, C_k\}$ and let $C_i$ be any one of such cycles. There are three possibilities:

\noindent Case 1. There is an edge between $C_i$ and another odd cycle say $C_j$. Let $u\in C_i$, $v\in C_j$ and $uv\in E(G)$. We put one of $u$ or $v$ in $D$. It is easy now to find $(|C_i|+|C_j|)/2$ suitable vertices from $C_i\cup C_j$ in order to activate whole vertices of $C_i\cup C_j$. We add these vertices to $D$.

\noindent Case 2. The cycle $C_i$ is only adjacent to even cycle(s). Let $C_i$ be adjacent to $C_l$ where $C_l$ is an even cycle and it is not already activated. Let $u\in C_l$ be a vertex adjacent to $C_i$. We add $u$ to $D$. Now we can activate the whole $C_i\cup C_l$ using $(|C_i|+|C_l|-1)/2$ vertices including the vertex $u$.

\noindent Case 3. The cycle $C_i$ is only adjacent to a vertex say $v$ of $F$. In this case we put $v$ in $D$. Note that by $(|C_i|+1)/2$ vertices we can activate whole $C_i\cup \{v\}$.

\noindent So far we have activated all the odd cycles of $\{C_1, C_2, \ldots, C_k\}$ and whole vertices of the forest $F$ except possibly some of its leaf vertices. Also some of the even cycles of the collection are activated during the above steps. The remaining even cycles can easily be activated using half of their vertices. We add these new vertices by which we activate even cycles in $D$ too. The construction process of the dynamo $D$ is finished. We only have to check the leaf vertices of $F$. Let $v$ be any vertex of degree one in a connected component $T$ of $F$. We already know that $v$ has an active neighbor in the tree $T$. Since the degree of $v$ in the whole graph is at least two then it should have another neighbor in $C_1\cup \ldots \cup C_k$. But the latter set is activated by $D$. Hence we obtain two active neighbors for the vertex $v$. Note that the set $D$ has at most $n/2$ vertices by the way of its construction. This completes the proof of part (i).

\noindent To prove part (ii), let $D$ be any dynamo of size $dyn(G[S])$. Using $D$ we can activate all vertices of $G[S]$. Let $T$ be any connected component of $G\setminus S$. Any vertex of $T$ has an active neighbor since $S$ is a domination set of $G$. Now since $T$ is connected one extra vertex from $T$ is enough to activate all vertices of $T$. We conclude that there exists a dynamo of size $dyn(G[S])+c$, where $c$ is the number of components of $G\setminus S$.

\noindent Part (iii) can be derived from part (ii) and applying the proof of part (i)
for $G[S]$.
\end{proof}

\noindent In the following theorem by determining the minimum dynamo of Generalized Petersen graphs we show that the lower bound of Corollary \ref{regular} is tight. We first recall the
definition of Generalized Petersen graphs $GP(n,k)$. For any positive integers $n$ and $k$ with $k\leq n-2$ and $gcd(n,k)=1$
we define $GP(n,k)$ as follows. The vertex set of $GP(n,k)$ is $\{a_1, a_2, \ldots, a_n\}\cup \{b_1, b_2, \ldots, b_n\}$. The edges of $GP(n,k)$ are $a_ib_i$ for $i=1, \ldots, n$, $a_1a_n$, $a_ja_{j+1}$ for $j=1, \ldots, n-1$ and also $a_ib_j$ for any $i$ and $j$ such that $|i-j|=k$. We note that $GP(n,k)$ is a cubic graph and its  order is $2n$. The following theorem determines the exact value of the minimum dynamic monopoly in $GP(n,k)$, where the threshold of any vertex is two.

\begin{thm}
The size of smallest dynamic monopoly in $GP(n,k)$ is $\lceil \frac{n+1}{2} \rceil$.\label{GP}
\end{thm}

\noindent \begin{proof}
Since $|GP(n,k)|=2n$, then by Corollary \ref{regular} it is enough to show that $GP(n,k)$ contains a dynamo
of size $\lceil \frac{2n+2}{4} \rceil = \lceil \frac{n+1}{2} \rceil$. Assume first that $n$ is an even number. In this case we use directly Theorem \ref{t=2}. The subset $S=\{a_1, a_2, \ldots, a_{n}\}$ forms a dominating set in $GP(n,k)$ and since $gcd(n,k)=1$ then $GP(n,k)\setminus S = G[b_1, \ldots, b_n]$ is a connected subgraph of the graph. From other side, $dyn(G[S])=n/2$. Hence by Theorem \ref{t=2} there exists a dynamo of size $\lceil \frac{n+1}{2} \rceil$ for $GP(n,k)$ when $n$ is even. In fact $\{a_1, a_3, \ldots, a_{n-1}, b_{n-1}\}$ is a dynamo of this size.

\noindent When $n$ is odd, it can be shown that $S'=\{a_1, a_3, a_5, \ldots, a_{n-2}, b_{n-1}\}$ is a dynamic monopoly in $GP(n,k)$. The argument uses Theorem \ref{t=2} but the point is that using $S'$ we first activate all vertices of $\{a_1, a_2, \ldots, a_n\}$ and then use the connectivity of $GP(n,k)\setminus \{a_1, a_2, \ldots, a_n\}$.
\end{proof}

\section{Dynamo-unbounded families of graphs}

\noindent Consider an election where people votes YES or NO to a certain candidate. Any individual (represented by a vertex $v$ in the underlying network) votes YES if at least $t(v)$ number of her friends have decided to vote YES. A dynamic monopoly $D$ for the underlying network of this election has the property that if the vote of the members of $D$ is YES then the whole community will eventually vote to that candidate. The following strategic question arises. If the population of the community increases then does it imply that the size of the smallest dynamo too increases (as a function of the size of community)? Another example where the same question becomes important is the adoption of a new product in viral marketing (for a formulation of viral marketing in terms of dynamic monopolies see \cite{DRi}). In the following by introducing the concept of dynamo-unbounded families we present a method to analyze the question we mentioned above.

\noindent By a threshold pattern we mean any threshold assignment $\tau$ such that for any graph $G$ and any $v\in V(G)$, $\tau$ assigns a non-negative value $t(v)$ such that if $\sigma$ is any automorphism of $G$ with $\sigma(u)=v$ for some vertex $u\in G$, then $t(v)=t(u)$. Without loss of generality we may restrict a threshold pattern $\tau$ so that $t(v)\leq d_G(v)$. The common examples are when $t(v)$ is a function of $d_G(v)$ for any vertex $v$. In this section by a family $\mathcal{F}$ we mean any set of graphs equipped with a threshold pattern. Such a family is called dynamo-unbounded if there exists a function $f(x)$ satisfying $f(x)\rightarrow \infty$ as $x\rightarrow \infty$ such that for any graph $G$ from $\mathcal{F}$ one has $f(n)\leq dyn(G)$, where $n=|G|$. Corollary \ref{regular} implies that the family of $2r+1$-regular graphs with threshold $t(v)=r+1$ for any vertex, is dynamo-unbounded family. In this section we obtain more results concerning dynamo-unbounded graphs.

\noindent In the following corollary we denote the edge density of a graph $G$ by $\epsilon(G)$
which is defined as $\epsilon(G)=|E(G)|/|V(G)|$.

\begin{cor}
Let $(G, \tau)$ be a graph of order $n$. Set $t=\min \{t(v):
v\in V(G)\}$. Then $$n(1-\frac{\epsilon(G)}{t})\leq dyn(G).$$\label{epsilon}
\end{cor}

\noindent \begin{proof}
Let $D$ be any dynamo of size $k=dyn(G)$ and $H=G\setminus K$. Then by part (ii) of Theorem \ref{bound1},
$(n-k)t\leq |E(G)|$. This easily implies the desired inequality.
\end{proof}

\noindent The following corollary follows immediately from Corollary \ref{epsilon}.

\begin{cor}
Let ${\mathcal{F}}$ be any family of graphs such that for some positive constant $\delta$, $\min\{t(v): v \in G\} \geq \epsilon(G) + \delta$ for any graph $G\in {\mathcal{F}}$. Then ${\mathcal{F}}$ is dynamo-unbounded.\label{cor3}
\end{cor}

\noindent The following theorem concerns graphs with probabilistic thresholds. In proving the following theorem we shall make use of the following concentration result of McDiarmid \cite{Mc}. Let $X_1, \ldots, X_n$ be a sequence of nonnegative independent random variables and set $X=\sum X_i$. Then for any $\lambda\geq 0$ $$\Bbb{P}(X\leq \Bbb{E}(X)-\lambda)\leq e^{-\frac{\lambda^2}{2\sum_{i=1}^n \Bbb{E}(X_i^2)}}.$$

\noindent \begin{thm}
Assume that any vertex $v$ of $G$ chooses a random threshold $a\leq i\leq b$ with probability $p_i$ where $p_i$ is independent of $v$. Set $\alpha=\sum_{a}^b ip_i$. If $\alpha > \epsilon(G)$ then for any positive constant $\delta$, with high probability no dynamo of $G$ contains less than $n^{1-\delta}$ elements.\label{homog}
\end{thm}

\noindent \begin{proof}
Let $n$ be a sufficiently large integer so that $\alpha /(\alpha -\epsilon) < n^{\delta}$. Let $D$ be any subset of vertices of cardinality $k$ where $k<n^{1-\delta}$ and set $H=G\setminus D$.
For any vertex $v\in H$ define a random variable $X_v$ as the threshold $t(v)$ chosen by the vertex $v$. Set $X=\sum_{v\in H} X_v$. We have the following by Theorem \ref{bound1}, part (ii).
$$\Bbb{P}(D~is~a~dynamo)\leq \Bbb{P}(X\leq |E(G)|).$$
\noindent For $X$ we have the following information

$\Bbb{E}(X)=(n-k)\sum_{i\in [a,b]} ip_i=\alpha(n-k)$~~~and~~~
$\sum_{v\in H} \Bbb{E}(X_v^2)=(n-k)\sum_{i\in [a,b]} i^2p_i$.

\noindent Also write $\beta = \sum_{i\in [a,b]} i^2p_i$ for simplicity. By our hypothesis $k<n^{1-\delta}$ which by $\alpha /(\alpha -\epsilon) < n^{\delta}$ imply $k<n(1-\epsilon/\alpha)$ or $\Bbb{E}(X) - |E(G)|>0$. Now by considering $\lambda = \Bbb{E}(X) - |E(G)|$ we may use the above-mentioned result of McDiarmid, since $\lambda >  0$. We obtain

$$\hspace{-2cm}\Bbb{P}(D~is~a~dynamo)\leq \Bbb{P}(X\leq |E(G)|)$$
$$\hspace{2cm}\leq \Bbb{P}(X\leq \Bbb{E}(X) - \lambda)$$
$$\hspace{1cm}\leq e^{-\frac{\lambda^2}{2\sum_{v\in H} \Bbb{E}(X_v^2)}}$$
$$\hspace{1.3cm}= e^{-\frac{(\alpha (n-k) - n\epsilon(G))^2}{2(n-k) \beta}}$$
$$\hspace{1cm}= e^{-\frac{(n(\alpha-\epsilon)-k\alpha)^2}{2(n-k) \beta}}$$
$$\hspace{1cm}\leq e^{-\frac{(n(\alpha-\epsilon)-k\alpha)^2}{2n\beta}}$$
$$\hspace{1.6cm}\leq e^{-\frac{(n(\alpha-\epsilon)-\alpha n^{1-\delta})^2}{2n\beta}}.$$

\noindent The latter inequality implies that for some positive constants $c'$ and $c$, $\Bbb{P}(D~is~a~dynamo)\leq e^{-\frac{c'n^2}{2n\beta}}=e^{-cn}$. Therefore the probability that there exists a monopoly of size $k$ is at most
$$n^ke^{-cn}=e^{k\ln n -cn}.$$

\noindent Now since $k\leq n^{1-\delta}$ then $e^{k\ln n -cn}=o(1)$. This shows that no subset of cardinality less than $n^{1-\delta}$ is a dynamo. This completes the proof.
\end{proof}

\noindent We call any graph satisfying the assumptions of Theorem \ref{homog}, a homogenious society. The proof of Theorem \ref{homog} shows that if $G$ is a homogenious society and $D$ any subset of $G$ with $|D|\leq n(1-\epsilon / \alpha)$ then with high probability $D$ is not a dynamo. We pose the following question: Is it true that with high probability $dyn(G)\geq n(1- \epsilon / \alpha)$ for any homogenious society?

\noindent A result related to the concept of dynamo-unbounded graphs is that of \cite{B}. In \cite{B} the reversible version of the model we studied in this paper has been considered. At each time step of the process any vertex updates its state as follows. Any vertex $v$ takes a new state which is the state of the majority of its neighbors. In case that the number of active neighbors is the same as the number of non-active neighbors of $v$ then the state of $v$ is remained unchanged. Berger proved that for any $n$ there exists a graph of more than $n$ vertices which contains a dynamic monopoly of at most 18 vertices.

\noindent The following proposition shows that the result of Corollary \ref{cor3} is the  best possible.

\begin{prop}
For any positive integers $r$ and $n$ with $r|n$, there exists a $2r$-regular graph on $n$ vertices which contains a dynamo of size $r$, where the threshold of any vertex is taken $r$.
\end{prop}

\noindent \begin{proof}
Write $n=rq$ for some $q>0$. Let also $C_1, \ldots, C_q$ be $q$ vertex disjoint copies of $\overline{K_r}$ where
$\overline{K_r}$ is the empty graph on $r$ vertices. Denote the vertex set of $C_i$ by $V_i$. We define a graph $G$ as follows. The vertex set of $G$ is $V_1\cup V_2 \cup \ldots \cup V_q$.
In $G$ the subgraph induced on $V_i\cup V_{i+1}$ for any $i=1, 2, \ldots, q$ (when $i=q$ we take $q+1$ as 1) is a complete bipartite subgraph whose
bipartition sets are $V_i$ and $V_{i+1}$. Set the threshold of each vertex of $G$ as $r$. It is easily seen that $dyn(G)=\epsilon(G)=r$.
\end{proof}

\section{Dynamic monopolies in line graphs}

\noindent By the line graph of a graph $G$ denoted by $L(G)$ we mean a graph
whose vertex set is the edge set of $G$ where two vertices $e$ and $e'$ of $L(G)$
(as two edges in $G$) are adjacent if and only if $e$ intersects $e'$ in $G$. A dynamic monopoly of $L(G)$ can be considered as the dynamic monopoly of the edges of $G$. In this section we study dynamic monopolies in line graphs when $G$ is a regular graph where our studies will be in terms of the edges of $G$ instead of working with vertices of $L(G)$. Note that if $e$ is any edge between two vertices $u$ and $v$ in a graph $G$ then the degree of $e$ as a vertex of $L(G)$ is $d_G(u)+d_G(v)-2$.

\begin{thm}
Let an $r$-regular graph $G$ with an assignment of thresholds to the edges of $G$ be given. Set $t=\min\{t(e):e\in E(G)\}$. Let $D\subseteq E(G)$ be a dynamic monopoly of size $k$ in $L(G)$. Then $$k\geq \lfloor \frac{4(t-r+1)n+(2r-t)^2}{8}\rfloor .$$\label{edgecomplete}
\end{thm}

\noindent \begin{proof}
\noindent Since $G$ is $r$-regular then $L(G)$ is $2r-2$-regular. Set $H=G\setminus D$. The graph $H$ has $n$ vertices and $|E(G)|-k=rn/2 -k$ edges. Since $D$ is a dynamo then there exists $e_1\in E(H)$ such that $d_H(e_1)=d_G(e_1)-d_D(e_1)\leq d_G(e_1)-t$. Note that if $e=uv$ then $d_H(e_1)=d_H(u)+d_H(v)-2$.

\noindent Set now $H_1=H\setminus \{u,v\}$ we have $|H_1|=n-2$ and when we remove $u$ and $v$ from $H$, we lose exactly $d_H(u)+d_H(v)-1$ edges from $H$.

$$|E(H_1)| = |E(H)|-(d_H(u)+d_H(v)-1)$$
$$\hspace{-0.25cm}=|E(H)|-d_H(e_1)-1$$
$$\hspace{0.4cm}\geq |E(H)| - d_G(e_1) + t -1$$
$$\hspace{0.8cm}\geq |E(H)| - (2r-2) + t-1.$$

\noindent We repeat the above technique and obtain $H_i$ on $n-2i$ vertices such that
$$|E(H_i)|\geq |E(H)| -i(2r-2)+i(t-1).$$
\noindent Now we use the obvious upper bound $|E(H_i)|\leq {n-2i \choose 2}$ and
obtain the following inequalities for any $i$ $$|E(H)| \leq i(2r-t-1)+\frac{4i^2-4in+2i+n^2-n}{2}$$
$$|E(H)| \leq 2i^2+i(2r-t-1-2n+1)+(n^2-n)/2.$$

\noindent The value in the right hand of the above inequality minimizes at
$i=\frac{2n-2r+t}{4}$. Its minimum value is $\frac{8rn-4nt-4n+4rt-4r^2-t^2}{8}$.
It turns out that
$$k\geq \lfloor \frac{4(t-r+1)n-4rt+4r^2+t^2}{8} \rfloor.$$
\end{proof}

\noindent The following theorem is concerning the line graphs of bipartite graphs.

\begin{thm}
Let $G$ be an $r$-regular bipartite graph on $n$ vertices and $t$ an assignment of thresholds to the edges of $G$. Set $t=\min\{t(e):e\in E(G)\}$. Let $D\subseteq E(G)$ be a dynamic monopoly of size $k$ in $L(G)$. Then $$k\geq \frac{n(2t-2r+2)+(2r-t)^2-4r+2t}{4} + \epsilon$$
\noindent where $\epsilon=1/4$ if $n-2r+t+1$ is an even integer and $\epsilon =0$ otherwise.\label{edgebip}
\end{thm}

\noindent \begin{proof}
The proof is similar to the one of Theorem \ref{edgecomplete}. Since $G$ is regular then each bipartition of $G$ contains $n/2$ vertices. There exists $e_i\in E(H_{i-1})$ such that
$d_{H_{i-1}}(e_i)\leq d_G(e_i)-t$. Set $H_i=H_{i-1}\setminus \{u_i,v_i\}$ where $e_i=u_iv_i$. We obtain $|H_i|=n-2i$ and
$$|E(H_i)|\geq |E(H)| -i(2r-2)+i(t-1).$$
\noindent We have now $|E(H)| \leq |E(H_i)|+i(2r-t-1)$.
From other side $|E(H_i)|\leq (\frac{n}{2} -i)^2$, since $H_i$ is a bipartite graph.

\noindent We have now $$|E(H)|\leq i^2+i(2r-t-n-1)+\frac{n^2}{4}.$$
\noindent The minimum value of the right hand term in the above inequality
is $\frac{n^2}{4}-\frac{(n-2r+t+1)^2}{4}$ and it is achieved when $i=(n+t-2r+1)/2$ is an integer, i.e. when $n+t-2r+1$ is even. Namely when $n+t-2r+1$ is even then $|E(H)|\leq \frac{n^2}{4}-\frac{(n-2r+t+1)^2}{4}$. But when $n+t-2r+1$ is odd then $|E(H)|\leq \frac{n^2}{4}-\frac{(n-2r+t+1)^2}{4}+1/4$. Therefore $|E(H)|\leq \frac{n^2}{4}-\frac{(n-2r+t+1)^2}{4}+ \varphi$, where $\varphi=0$ when $n+t-2r+1$ is even
and $\varphi=1/4$ when $n+t-2r+1$ is odd. Also $|E(H)|=|E(G)|-k=rn/2 - k$. We have the following
$$k\geq rn/2 - |E(H)| \geq \frac{rn}{2} - \frac{n^2}{4} + \frac{n^2+4r^2+t^2+1-4nr+2nt+2n-4rt-4r+2t}{4}-\varphi$$
$$\hspace{.9cm}\geq \frac{n(2t-2r+2)+(4r^2+t^2-4rt-4r+2t+1)}{4}-\varphi$$
$$\hspace{-.10cm}= \frac{n(2t-2r+2)+(2r-t)^2-4r+2t}{4}+\frac{1}{4}-\varphi.$$
\noindent By taking $\epsilon=1/4-\varphi$ the proof completes.
\end{proof}\\

\noindent The following result deals with the line graphs of complete graph and complete bipartite graph $K_{n,n}$, with constant edge thresholds $n-2$ (for $K_n$) and $n-1$ (for $K_{n,n}$). Note that $L(K_{n,n})=K_n\Box K_n$, where $\Box$ denotes the Cartesian product of graphs.

\begin{cor}

(i) Any dynamic monopoly for the edges of $K_n$ has at least $\lfloor n^2/8 \rfloor$ vertices.

(ii) The size of smallest dynamo in $K_n\Box K_n$ is $\lfloor \frac{n^2}{4} \rfloor$.
\end{cor}

\begin{figure}[ht]
\[
{\large \textcolor{red}{\begin{tabular}
{|p{.25cm}|p{.25cm}|p{.25cm}|
p{.25cm}|p{.25cm}|p{.25cm}|p{.25cm}|}
\hline&&&&{\bf *}&{\bf *}&{\bf *}\\[-0.1eM]
\hline&&&&&{\bf *}&{\bf *}\\[-0.1eM]
\hline&&&&&&{\bf *}\\[-0.1eM]
\hline&&&&&&\\[-0.1eM]
\hline{\bf *}&&&&&&\\[-0.1eM]
\hline{\bf *}&{\bf *}&&&&&\\[-0.1eM]
\hline{\bf *}&{\bf *}&{\bf *}&&&&\\[-0.0eM]
\hline
\end{tabular}}}\vspace{0.6cm}
\hspace{2cm}
{\large \textcolor{red}{\begin{tabular}
{|p{.25cm}|p{.25cm}|p{.25cm}|p{.25cm}|
p{.25cm}|p{.25cm}|p{.25cm}|p{.25cm}|}
\hline&&&&&{\bf *}&{\bf *}&{\bf *}\\[-0.1eM]
\hline&&&&&&{\bf *}&{\bf *}\\[-0.1eM]
\hline&&&&&&&{\bf *}\\[-0.1eM]
\hline&&&&&&&\\[-0.1eM]
\hline&&&&{\bf *}&&&\\[-0.1eM]
\hline{\bf *}&&&&&{\bf *}&&\\[-0.1eM]
\hline{\bf *}&{\bf *}&&&&&{\bf *}&\\[-0.1eM]
\hline{\bf *}&{\bf *}&{\bf *}&&&&&{\bf *}\\[-0.0eM]
\hline
\end{tabular}}}\vspace{0.6cm}
%\caption{A monopoly of size 16 in $K_8\Box K_8$}
\label{latin}
\]\caption{Minimum dynamos for $K_7\Box K_7$ and $K_8\Box K_8$}
\end{figure}

\noindent \begin{proof}
\noindent By applying Theorem 4 for $G=K_n$, $r=n-1$ and $t=n-2$, we obtain $dyn(K_n)\geq \lfloor n^2/8 \rfloor$. To prove (ii) we note that a dynamic monopoly in $K_n\Box K_n$ is equivalent to an edge dynamic monopoly in $K_{n,n}$. In this case $r=n$ and $t=n-1$ and the lower bound in Theorem \ref{edgebip} shows $k\geq \lfloor n^2/4 \rfloor$. In the following we obtain a dynamo of size $\lfloor \frac{n^2}{4} \rfloor$ for $K_n\Box K_n$. Consider the vertex set of $K_n\Box K_n$ as an $n\times n$ square array, where each vertex is identified by a position say $(i,j)$ in the array (i.e. the $i$-th row and $j$-th column). First let $n$ be an odd integer and write $n=2k+1$. In this case our dynamo $D$ consists of two triangular subarrays (see Figure 1 for $n=7$) in down-left and top-right parts of the whole array. The height and side of these two triangular arrays are $k$. It can be easily checked that the resulting subarray is a dynamo indeed and it has $k(k+1)=(n^2-1)/4$ entries. In fact the positions $(1,1)$ and $(n,n)$ are the first vertices which become active. Then the first row and column and the last row and column of the array become active. We reach at an array of size $(n-2)\times (n-2)$ where in addition to $D$, all the vertices in the first and last rows and columns are also activated. The rest of the array becomes active inductively in a similar manner. This proves the theorem for odd $n$.

\noindent When $n$ is even of the form $n=2k$ we consider two similar triangular subarrays except that the height and side of these triangles are $k-1$ (see Figure 1 for $n=8$, where $k=4$). In addition to the vertices of these subarrays we also consider a subset of vertices on the principal diagonal of the array in our dynamo consisting of the positions $(k+1, k+1), (k+2, k+2), \ldots, (n, n)$. The argument that the chosen vertices form a dynamo is similar to the previous one and we omit its proof.
\end{proof}

%%%%%%%%%%%%%%%%%%%%%%%%%%%%%%%%%%%%%%%%%%%%%%%%%%%%%%%%%%%%%%%%%%%%%%%%%%%%%%%%%%%

\end{document}